\documentclass[11pt]{amsart}
\usepackage{a4wide}
\usepackage[utf8]{inputenc}
\usepackage{amsmath}
\usepackage{amssymb}
\usepackage{amsfonts}
\usepackage{amsopn}
\usepackage{graphicx}
\usepackage{enumerate}
\usepackage{color}
\usepackage{mathtools}
\usepackage{mathrsfs}
\usepackage{bbm}
\usepackage{yhmath}
\usepackage[colorlinks,linkcolor=blue,anchorcolor=blue,citecolor=blue]{hyperref}
\usepackage[all,pdf]{xy} 
\usepackage{tikz}

\newtheorem{theorem}{Theorem}[section]
\newtheorem{proposition}[theorem]{Proposition}
\newtheorem{lemma}[theorem]{Lemma}

\newtheorem{corollary}[theorem]{Corollary}
\theoremstyle{remark}
\newtheorem{remark}[theorem]{Remark}
\newtheorem{example}[theorem]{Example}

\numberwithin{equation}{section}

\newcommand{\mrm}{\mathrm}

\newcommand{\mscr}{\mathscr}
\newcommand{\set}[1]{\left\{ #1 \right\}}

\newcommand{\R}{\mathbb{R}}
\newcommand{\Q}{\mathbb{Q}}

\newcommand{\N}{\mathbb{N}}

\newcommand{\f}{\infty}
\newcommand{\la}{\lambda}
\newcommand{\al}{\alpha}
\newcommand{\ga}{\gamma}

\newcommand{\de}{\delta}
\newcommand{\sse}{\subset}
\newcommand{\sm}{\setminus}

\newcommand{\La}{\Lambda}
\newcommand{\Om}{\Omega}
\newcommand{\om}{\omega}
\newcommand{\V}{\mscr{V}}
\newcommand{\W}{\mscr{W}}
\newcommand{\lle}{\preccurlyeq}
\newcommand{\lge}{\succcurlyeq}

\allowdisplaybreaks

\title[Self-similar sets with common points]{On a class of self-similar sets which contain finitely many common points}

\author{Kan Jiang}
\address[K. Jiang]{School of Mathematics and Statistics, Ningbo University, Ningbo 315211, People's Republic of China}
\email{jiangkan@nbu.edu.cn}

\author{Derong Kong}
\address[D. Kong]{College of Mathematics and Statistics, Center of Mathematics, Chongqing University, Chongqing 401331, People's Republic of China}
\email{derongkong@126.com}

\author{Wenxia Li}
\address[W. Li]{School of Mathematical Sciences, Key Laboratory of MEA (Ministry of Education) \& Shanghai Key Laboratory of PMMP,  East China Normal University, Shanghai 200241, People's Republic of China}
\email{wxli@math.ecnu.edu.cn}

\author{Zhiqiang Wang}
\address[Z. Wang]{School of Mathematical Sciences, Key Laboratory of MEA (Ministry of Education) \& Shanghai Key Laboratory of PMMP, East China Normal University, Shanghai 200241, People's Republic of China}
\email{zhiqiangwzy@163.com}

\begin{document}

\begin{abstract}
For $\lambda\in(0,1/2]$ let $K_\lambda \subset\R$ be a self-similar set generated by the iterated function system $\{\lambda x, \lambda x+1-\lambda\}$.
Given $x\in(0,1/2)$, let $\Lambda(x)$ be the set of $\lambda\in(0,1/2]$ such that $x\in K_\lambda$.
In this paper we show that $\Lambda(x)$ is a topological Cantor set having zero Lebesgue measure and full Hausdorff dimension. Furthermore, we show that for any $y_1,\ldots, y_p\in(0,1/2)$ there exists a full Hausdorff dimensional set of $\lambda\in(0,1/2]$ such that $y_1,\ldots, y_p \in K_\lambda$.
\end{abstract}

\keywords{Hausdorff dimension, thickness, self-similar set, Cantor set, intersection.}

\subjclass[2020]{Primary: 28A78; Secondary: 28A80, 37B10, 11A63 }

\maketitle

\section{Introduction}

For $\lambda\in(0,1/2]$ let $K_\lambda$ be the \emph{self-similar set} generated by the \emph{iterated function system} (simply called, \emph{IFS})
$\set{f_{\lambda,d}(x)=\lambda x + d(1-\lambda): d=0,1 }.$
Then $K_\lambda$ is the unique nonempty compact set satisfying (cf.~\cite{Hutchinson_1981})
\begin{equation}\label{eq:self-similar-set}
 K_\lambda= f_{\la,0}(K_\la)\cup f_{\la,1}(K_\la)=\set{ (1-\la)\sum_{n=1}^\f i_n\la^{n-1}: i_n\in\set{0,1}~\forall n\ge 1 }.
\end{equation}
Clearly, $1-\lambda$ is chosen so that the convex hull of $K_\la$ is the unit interval $[0,1]$ for all $\la\in(0,1/2]$. Then $0$ and $1$ are common points of  $K_\la$ for all $\la\in(0,1/2]$.
For other $x\in(0,1)$ it is natural to ask how likely  the self-similar sets $K_\la, \la\in(0,1/2]$ contain  the  common point $x$? Or  even ask how likely the self-similar sets $K_\la, \la\in(0,1/2]$ contain any given points $y_1,\ldots, y_p\in(0,1)$?
These questions are motivated by the work of Boes, Darst and Erd\H{o}s \cite{BDE-1981}, in which they considered a   class of fat Cantor sets $C_\la$ with positive Lebesgue measure. They showed that for a given point $x\in(0,1)$   the set of parameters $\la\in(0,1/2)$ such that $x\in C_\la$ is of first category.

Given $x\in[0,1]$, let
\begin{equation}\label{eq:definition-La}
 \Lambda(x) := \set{ \lambda \in (0,1/2]: x \in K_\lambda }.
\end{equation}
Then $\La(x)$ consists of all  $\la\in(0,1/2]$ such that $x$ is the common point of  $K_\la$.
Note that  $K_\la$ is symmetric, i.e., $x\in K_\la$ if and only if $1-x\in K_\la$. Then
$\La(x)=\La(1-x)$ for any $x\in [0,1].$
So, we only need to consider $x\in[0,1/2]$. Note that $\La(0)=(0,1/2]$, and $\La(1/2)=\set{1/2}$. So, it is interesting  to study $\La(x)$ for $x\in(0,1/2)$.

Recall that a set $F\subset\R$ is called a \emph{Cantor set} if it is a non-empty compact set containing neither interior nor isolated points. Our first result considers the topology of $\La(x)$.
\begin{theorem}\label{main:topology}
  For any $x\in (0,1/2)$  the  set $\Lambda(x)$ is a Cantor set with $\min \Lambda(x) = x$ and $\max \Lambda(x) = 1/2$.
\end{theorem}

By Theorem  \ref{main:topology} it follows that $\La(x)$ is a fractal  set for any $x\in(0,1/2)$.  {Our next result considers the Lebesgue measure and fractal dimension of $\La(x)$.}

\begin{theorem}\label{main:measure-dim}
  For any $x \in (0, 1/2)$ the set $\La(x)$ is a Lebesgue null set of full Hausdorff dimension.
 Furthermore,
\[ \lim_{\de\to 0^+}\dim_H(\La(x)\cap(\la-\de,\la+\de))=\frac{\log 2}{-\log \la}\quad\forall~\la\in\La(x), \]
where $\dim_H$ denotes the Hausdorff dimension.
\end{theorem}
In 1984, Mahler \cite{Mahler-1984} proposed the problem on studying how well elements in the middle third Cantor set $K_{1/3}$  can be approximated by rational numbers in it, and by rational numbers outside of it. Some recent progress on this problem can be found in   \cite{Levesley-Salp-Velani-2007, Fishman-Simmons-2014, Fishman-Simmons-2015, Schleischitz-2021} and the references therein. On the other hand, this question also motivates the study of rational numbers in a fractal set (cf.~\cite{Shparlinski-2021, Wall-1990}). As a corollary of Theorem \ref{main:measure-dim} we show that for Lebesgue almost every $\lambda\in(0,1/2)$ the Cantor set $K_\lambda$ contains only two rational numbers $0$ and $1$.

\begin{corollary}\label{cor:irrational-numbers}
  For Lebesgue almost every $\lambda \in (0,1/2]$  the set $K_\lambda\setminus\set{0,1}$ contains only irrational numbers.
\end{corollary}
\begin{proof}
  By Theorem \ref{main:measure-dim} it follows that $\La(x)$ has zero Lebesgue measure for any $x\in(0,1)$.
  But if $K_\lambda \setminus \{0,1\}$ contains a rational number, then $\la \in \bigcup_{x\in\Q\cap(0,1)}\La(x)$ which has zero Lebesgue measure.
\end{proof}

Given  $y_1,\ldots, y_p\in(0,1/2)$, by Theorems \ref{main:topology} and \ref{main:measure-dim} it follows that the intersection $\bigcap_{i=1}^p\La(y_i)$ is small from the topological and Lebesgue measure perspectives.  {On the other hand}, by using the \emph{thickness} method introduced by Newhouse \cite{Newhouse-1970} we {can} show that  $\bigcap_{i=1}^p\La(y_i)$ contains a sequence of Cantor sets whose thickness can be arbitrarily large, {and from this we conclude that the intersection  $\bigcap_{i=1}^p\La(y_i)$ has   full Hausdorff dimension.}
\begin{theorem}\label{main:intersection}
  For any  points $y_1, y_2, \ldots, y_p \in (0, 1/2)$  we have
  \[\dim_H \bigcap_{i=1}^p\Lambda(y_i) =1.\]
\end{theorem}

Recently, the first three authors studied in \cite{Jiang-Kong-Li-2021} analogous objects but with  different family of self-similar sets (their self-similar sets have different convex hulls).
Furthermore, Theorem \ref{main:intersection} shows that the intersection of any  finitely many $\La(y_i)$ has full Hausdorff dimension, while in \cite[Theorem 1.5]{Jiang-Kong-Li-2021}  their method can only  prove {this result for  the intersection of two associated  sets.}

The rest of the paper is organized as follows. In Section \ref{sec:topology} we prove Theorem \ref{main:topology} for the topology of $\La(x)$; and in Section \ref{sec:measure-dim} we investigate the local Hausdorff dimension of $\La(x)$ and prove Theorem \ref{main:measure-dim}. In Section \ref{sec:intersection} we consider the intersection $\bigcap_{i=1}^p\La(y_i)$ and prove Theorem \ref{main:intersection}; and in the final section we  make some further remarks.

\section{Topological properties of $\La(x)$}\label{sec:topology}
In this section we investigate the topology of $\La(x)$, and prove Theorem \ref{main:topology}.
{Given $x\in(0,1/2)$, note by (\ref{eq:self-similar-set}) that for each $\la\in\La(x)\setminus\set{1/2}$ there exists a   unique  sequence $(d_i)\in\set{0,1}^\N$ such that $x=(1-\la)\sum_{n=1}^\f d_n\la^{n-1}$.  We will show that $\La(x)$ is homeomorphic to a subset in the symbolic space $\set{0,1}^\N$, and then the topological properties of $\La(x)$ can be   deduced by studying the corresponding symbolic set.}

First we recall some terminology from symbolic dynamics (cf.~\cite{Lind_Marcus_1995}). Let $\set{0,1}^\N$ be the set of all infinite sequences of zeros and ones.  For a \emph{word} we mean a finite string of zeros and ones. Let $\set{0, 1}^*$ be the set of all words over the alphabet $\set{0,1}$ together with the empty word $\epsilon$. For two words $\mathbf c=c_1\ldots c_m, \mathbf d=d_1\ldots d_n$ from $\set{0,1}^*$ we write $\mathbf{cd}=c_1\ldots c_md_1\ldots d_n$ for their concatenation.
In particular, for $n\in\N$ we denote by $\mathbf c^n$ the $n$-fold concatenation of $\mathbf c$ with itself, and by $\mathbf c^\f$ the periodic sequence with {period} block $\mathbf c$. Throughout the paper we will use lexicographical order `$\prec, \lle, \succ$' or `$\lge$' between sequences and words. For example, for two sequences $(c_i), (d_i)\in\set{0,1}^\N$, we say $(c_i)\prec (d_i)$ if $c_1<d_1$, or there exists $n\in\N$ such that $c_1\ldots c_n=d_1\ldots d_n$ and $c_{n+1}<d_{n+1}$. For two words $\mathbf c, \mathbf d$, we say $\mathbf c\prec \mathbf d$ if $\mathbf c 0^\f\prec \mathbf d 0^\f$.

Let $\la\in(0,1/2]$. We define the coding map
$\pi_\lambda: \set{0,1}^\N\to K_\lambda$ by
\begin{equation}\label{eq:projection-map}
  \pi_\lambda((i_n))  = \lim_{n \to \f} f_{\lambda, i_1} \circ f_{\lambda,i_2} \circ \cdots \circ f_{\lambda, i_n} (0)   = (1-\lambda) \sum_{n=1}^{\f} i_n \lambda^{n-1}.
\end{equation}
If $\la\in(0,1/2)$, then the IFS $\set{f_{\la,d}(x)=\la x+d(1-\la): d=0,1}$ satisfies the strong separation condition, and thus the map $\pi_\la$
is bijective. If $\la=1/2$, then  $\pi_{1/2}$ is bijective up to a countable set. The map $\pi_\la$ defined in (\ref{eq:projection-map}) naturally induces a function with two parameters:
\begin{equation}\label{eq:projection-two}
\Pi: \set{0,1}^\N\times (0,1/2]~\to~[0,1];\quad ((i_n), \la)\mapsto \pi_\la((i_n)).
\end{equation}
Note that the symbolic space $\{0,1\}^\N$ becomes a  compact metric space under  the metric
\begin{equation}\label{eq:metric-rho}
\rho((i_n), (j_n))=2^{-\inf\set{n\ge 1: i_n\ne j_n}}.
\end{equation}
Equipped with the product topology on $\set{0,1}^\N\times (0,1/2]$ we show that $\Pi$ is continuous.

\begin{lemma}\label{lem:property-Pi}
The function $\Pi$ is continuous. Furthermore,
\begin{itemize}
\item [{\rm(i)}]   for $\lambda \in (0,1/2]$  the   function  $\Pi(\cdot,\la)$ is   increasing with respect to the lexicographical order, and is strictly increasing if $\la\in(0,1/2)$;

\item [{\rm(ii)}]  if $0^\f\prec (i_n)\lle 01^\f$, then $\Pi((i_n),\cdot)$ has positive derivative {in $(0,1/2)$}.
\end{itemize}
\end{lemma}

\begin{proof}
First we prove the continuity of $\Pi$. For any two points  $((i_n), \la_1), ((j_n), \la_2)\in\set{0,1}^\N\times(0,1/2]$ we have
\begin{equation}\label{eq:continuity-1}
  |\Pi((j_n), \la_2)-\Pi((i_n), \la_1)|\le  |\Pi((j_n), \la_2)-\Pi((i_n), \la_2)|+ |\Pi((i_n), \la_2)-\Pi((i_n), \la_1)|.
\end{equation}
Note that  if $\rho((j_n), (i_n))\le 2^{-m}$, then  $|\Pi((j_n), \la_2)-\Pi((i_n), \la_2)|\le \la_2^{m-1}\le 2^{1-m}$. So the first term in (\ref{eq:continuity-1}) converges to zero as  $\rho((j_n), (i_n))\to 0$.
Moreover, since the series $\Pi((i_n), \la)=(1-\la)\sum_{n=1}^\f i_n\la^{n-1}$ with parameter $\la$ converges uniformly in $(0,1/2]$, the second term in (\ref{eq:continuity-1}) also converges to zero as $|\la_2-\la_1|\to 0$. Therefore, $\Pi$ is continuous.

For (i) let $\la\in(0,1/2]$ and take two sequences  $(i_n), (j_n)\in\set{0,1}^\N$. Suppose
$(i_n)\prec (j_n)$. Then there exists $m\in\N$ such that ${ i_1\ldots i_{m-1} =j_1\ldots j_{m-1}}$ and $i_m<j_m$. This implies that
  \begin{align*}
    \Pi((i_n), \la) = (1-\lambda)\sum_{n=1}^{\f} i_n \lambda^{n-1}
     & \le (1-\lambda) \left( \sum_{n=1}^{m} i_n \lambda^{n-1} + \sum_{n=m+1}^{\f} \lambda^{n-1} \right) \\
    & \le (1-\lambda) \sum_{n=1}^{m} j_n \lambda^{n-1} \\
    & \le \Pi((j_n), \la),
  \end{align*}
  where the second inequality follows from $\la/(1-\la) \le 1$ for $\la\in(0,1/2]$, and this inequality is strict if $\la\in(0,1/2)$.

  For (ii)  let $(i_n)\in\set{0,1}^\N$ with $0^\f\prec (i_n)\lle 01^\f$. Then {$i_1=0$. So} for any $\la\in(0,1/2)$ we have $\Pi((i_n),\la)= (1-\lambda) \sum_{n=2}^{\f} i_n \lambda^{n-1}.$ This implies that
\begin{equation}\label{deriv}
\frac{\mrm{d} \Pi((i_n), \la)}{\mrm{d} \lambda} = \sum_{n=2}^{\f} n \left( \frac{n-1}{n} -\lambda \right) i_n \lambda^{n-2}>0,
\end{equation}
where the {strict} inequality follows since {$\la< 1/2$} and $(i_n)\succ 0^\f$. This completes the proof.
\end{proof}

Note that the map $\Pi$ defined in (\ref{eq:projection-two}) is surjective but not injective. Given $x\in[0,1]$, for $\la\in(0,1/2]$ we consider the horizontal fiber
\[
\Gamma_x(\la):=\Pi^{-1}(x)\cap\left(\set{0,1}^\N\times \set{\la}\right)=\set{((i_n), \la)\in\set{0,1}^\N\times(0,1/2]: (1-\la)\sum_{n=1}^\f i_n\la^{n-1}=x}.
\]
Then $\Gamma_x(\la)\ne\emptyset$ if and only if $\la\in\La(x)$.
Furthermore, by Lemma \ref{lem:property-Pi} (i) it follows that for any $\la\in\La(x)\cap(0,1/2)$ the fiber set $\Gamma_x(\la)$ consists of only one sequence; and  if $\la=1/2\in \La(x)$ then the set $\Gamma_x(1/2)$ {consists of}  at most two sequences. This defines a map
\[
\Psi_x: \La(x)\to \set{0,1}^\N;\quad \la\mapsto \Psi_x(\la),
\]
where $\Psi_x(\la)$ denotes the lexicographically largest sequence in $\Gamma_x(\la)$.
The sequence $\Psi_x(\la)$ is also called the \emph{greedy coding} of $x$ in base $\la$.

Given $x\in(0,1/2)$, we reserve the notation $(x_n):=\Psi_x(1/2)$ for the greedy coding of $x$ in base $1/2$. Then $(x_n)$ {begins with $0$ and} does not end with $1^\f$.
\begin{lemma}
  \label{lem:property-Psi-x}
  For any $x\in(0,1/2)$  the map $\Psi_x: \La(x)\to \Om(x)$ is a decreasing homeomorphism,  where
  \[
  \Om(x):=\set{(i_n)\in\set{0,1}^\N: (x_n)\lle (i_n)\lle 01^\f}.
  \]
\end{lemma}
\begin{proof}
Let $x\in(0,1/2)$.  By Lemma \ref{lem:property-Pi}  it follows that $\Psi_x$ is strictly decreasing. Observe that $x\notin K_\la$ for any $\la<x$. Then $\La(x)\subset[x, 1/2]$. Note that $\Psi_x(x)=01^\f$ and $\Psi_x(1/2)=(x_n)$. Since $\Psi_x$ is {monotonically} decreasing, we have
\begin{equation*}
(x_n)\lle \Psi_x(\la)\lle 01^\f\quad\forall ~\la\in \La(x).
\end{equation*}
So, $\Psi_x(\La(x))\subset\Om(x)$.

Next we show that $\Psi_x(\La(x))=\Om(x)$. Let $(i_n)\in\Om(x)$. Then  by Lemma \ref{lem:property-Pi} it follows that
\[
\Pi\left((i_n), \frac{1}{2}\right)\ge \Pi\left((x_n), \frac{1}{2}\right)=x \quad\textrm{and}\quad \Pi((i_n), \la)\searrow 0<x\quad\textrm{as }\la\searrow 0.
\]
So, by the continuity of $\Pi$ in Lemma \ref{lem:property-Pi} there must exist $\la\in(0,1/2]$ such that
\begin{equation}\label{eq:inclusion-1}
\Pi((i_n), \la)=x.
\end{equation}
If $\la\in(0,1/2)$, then (\ref{eq:inclusion-1}) gives that $\Psi_x(\la)=(i_n)$. If $\la=1/2$, then by (\ref{eq:inclusion-1}) and using $(i_n)\lge (x_n)$ we still have  $\Psi_x(\la)=(i_n)$.  This proves  $\Psi_x(\La(x))=\Om(x)$. Hence, $\Psi_x: \La(x)\to \Om(x)$ is a decreasing bijection.

To completes the proof it remains to prove the continuity of $\Psi_x$ and its inverse $\Psi_x^{-1}$. Since the proof for the continuity of $\Psi_x^{-1}$ is similar, we only prove it for $\Psi_x$. Take $\la_*\in\La(x)$.  Suppose $\Psi_x$ is not continuous at $\la_*$. Then there exists $N\in\N$ such that for any $\de>0$ we can find $\la\in\La(x)\cap(\la_*-\de, \la_*+\de)$ such that $|\Psi_x(\la)-\Psi_x(\la_*)|\ge 2^{-N}$. Letting $\de=1/k$ with $k=1,2,\ldots,$ we can find a sequence $(\la_k)\subset\La(x)$ such that
\begin{equation}
  \label{eq:continuity-2}
  \lim_{k\to\f}\la_k=\la_*\quad\textrm{and}\quad |\Psi_x(\la_k)-\Psi_x(\la_*)|\ge 2^{-N}\quad\forall ~k\ge 1.
\end{equation}
Write $\Psi_x(\la_k)=(i_n^{(k)})$ and $\Psi_x(\la_*)=(i_n^*)$. Then by (\ref{eq:continuity-2}) we have $i^{(k)}_1\ldots i^{(k)}_N\ne i_1^*\ldots i_N^*$ for all $k\ge 1$. Note that $(\set{0,1}^\N,\rho)$ is a compact metric space, where $\rho$ is defined in (\ref{eq:metric-rho}). So we can find a subsequence $\{k_j\}\subset\N$ such that the limit $\lim_{j\to\f} (i_n^{(k_j)})$ exists, say $(i_n')$. Then $i_1'\ldots i_N'\ne i_1^*\ldots i_N^*$. Observe that
\begin{equation}\label{eq:continuity-3}
\Pi((i_n^{(k_j)}), \la_{k_j})=x=\Pi((i_n^*), \la_*)\quad\forall j\ge 1.
\end{equation}
Letting $j\to\f$ in (\ref{eq:continuity-3}), by (\ref{eq:continuity-2}) and Lemma \ref{lem:property-Pi} it follows that
\begin{equation}\label{eq:continuity-4}
\Pi((i_n'), \la_*)=x=\Pi((i_n^*), \la_*).
\end{equation}
If $\la_*=1/2$, then $\la_{k_j}\le \la_*$ for all $j\ge 1$. Since $\Psi_x$ is decreasing,
it follows that $\Psi_x(\la_{k_j})\lge \Psi_x(\la_*)=(i_n^*)$ for all $j\ge 1$, and thus $(i_n')\lge (i_n^*)$. Note that $(i_n^*)$ is the greedy coding of $x$ in base $\la_*$. Then by (\ref{eq:continuity-4}) it follows that $(i_n')=(i_n^*)$, leading to a contradiction with  $i_1'\ldots i_N'\ne i_1^*\ldots i_N^*$. If $\la_*<1/2$, then (\ref{eq:continuity-4}) gives  that $(i_n')=(i_n^*)$. This again leads to a contradiction.
Therefore, $\Psi_x$ is continuous at $\la_*$. Since $\la_*\in\La(x)$ is  arbitrary, $\Psi_x$ is continuous in $\La(x)$.
This completes the proof.
\end{proof}

\begin{proof}[Proof of Theorem \ref{main:topology}]
Let $x\in(0,1/2)$. By Lemma \ref{lem:property-Psi-x} it follows that $\min\La(x)=\Psi_x^{-1}(01^\f)=x$ and $\max\La(x)=\Psi^{-1}_x((x_n))=1/2$. Observe that $\Om(x)$ is a Cantor set under the metric $\rho$ defined in (\ref{eq:metric-rho}), which means that $\Om(x)$ is a non-empty compact, perfect and totally disconnected set under  $\rho$.
Then by Lemma \ref{lem:property-Psi-x} we conclude that $\La(x)$ is also a Cantor set.
\end{proof}

\section{Lebesgue measure and Hausdorff dimension of $\La(x)$}\label{sec:measure-dim}

In this section we will prove Theorem \ref{main:measure-dim}, which states that for any $x\in(0,1/2)$ the set $\La(x)$ is a Lebesgue null set of full Hausdorff dimension.
The key ingredient in our proof of Theorem \ref{main:measure-dim} is Proposition \ref{th:local-dimension} (see below), which indicates that the  local Hausdorff dimension of $\La(x)$ at some $\la_0\in\La(x)$ is equal to the Hausdorff dimension of the self-similar set $K_{\la_0}$. This property on the interplay between the
‘parameter space’ (in this case, $\La(x)$) and the ‘dynamical space’ (in our case $K_\la$)  was
first observed by Douady \cite{Douady-1993} in the context of dynamics of real quadratic polynomials. A
similar result was proved by Tiozzo \cite{Tiozzo-2015}, who considers for $c\in\R$ the set of angles of external
rays which ‘land’ on the real slice of the Mandelbrot set to the right of $c$ (parameter space) and
the set of external angles which land on the real slice of the Julia set of the map $z\mapsto z^2+c$
(dynamical space), showing that these two sets have the same Hausdorff dimension. Some other similar results in different settings can be found in  \cite{Urbanski-1986, Bonanno-Carminati-Isola-Tiozzo-2013, Kong-Li-Lu-Wang-Xu-2020, Carminati-Tiozzo-2022}.

\begin{proposition}\label{th:local-dimension}
Let $x \in (0,1/2)$. Then  for any $\lambda \in \Lambda(x)$  we have
  \begin{equation}\label{eq:local-dimension}
   \lim_{\delta \to 0^+} \dim_H (\Lambda(x) \cap (\lambda-\delta,\lambda+\delta)) = \dim_H K_\lambda = -\frac{\log 2}{\log \lambda}.
  \end{equation}
\end{proposition}
The second equality in (\ref{eq:local-dimension}) is obvious, since for any $\la\in\La(x)$ the self-similar set $K_\la$ is generated by the IFS $\set{\la x, \la x+(1-\la)}$ satisfying the open set condition (cf.~\cite{Hutchinson_1981}). So it suffices to prove the first equality in (\ref{eq:local-dimension}).

\begin{lemma}
  \label{lem:local-dim-upper-bound}
   Let $x\in (0,1/2)$. Then for any $\lambda \in (x,1/2)$ we have
    \[ \dim_H (\Lambda(x) \cap [x,\lambda]) \le \dim_H K_\lambda. \]
\end{lemma}
\begin{proof}
  Let $\la\in(x, 1/2)$. Note by Lemma \ref{lem:property-Psi-x} that $\pi_\la\circ\Psi_x: \La(x)\cap[x,\la]\to K_\la$ is injective.
  By \cite[Proposition 3.3]{Falconer_1990} we only need to prove that the inverse map $(\pi_\la\circ\Psi_x)^{-1}$ is Lipschitz. In other words, it suffices to prove that for any $\la_1, \la_2\in\La(x)\cap[x, \la]$ we have
  \begin{equation}\label{eq:lipschitz-1}
   |\pi_\lambda(\Psi_x(\lambda_1))- \pi_\lambda(\Psi_x(\lambda_2))| \ge C |\lambda_1 -\lambda_2|,
   \end{equation}
  where $C>0$ is a constant {independent of $\la_1$ and $\la_2$}.

  Take $\lambda_1,\lambda_2 \in \Lambda(x)\cap [x,\lambda]$ with $\lambda_1 < \lambda_2$, and write $\Psi_x(\lambda_1) = (i_n)$, $\Psi_x(\lambda_2) = (j_n).$
  By Lemma \ref{lem:property-Psi-x} we have $i_1 =j_1=0$ and $(i_n) \succ (j_n)$. Then there exists $m \ge 2$ such that $i_1\ldots i_{m-1} = j_1\ldots j_{m-1}$ and $ i_m > j_m$.
  Note that \[ (1-\lambda_1) \sum_{n=2}^{\f} i_n \lambda_1^{n-1} =x =  (1-\lambda_2) \sum_{n=2}^{\f} j_n \lambda_2^{n-1}. \]
  Then
  \begin{align*}
    \frac{x(1-\lambda_1 -\lambda_2)}{\lambda_1 \lambda_2 (1-\lambda_1)(1-\lambda_2)}(\lambda_2 - \lambda_1) & = \frac{x}{\lambda_1(1-\lambda_1)} - \frac{x}{\lambda_2(1-\lambda_2)} \\
    & = \sum_{n=2}^{\f} i_n \lambda_1^{n-2} -\sum_{n=2}^{\f} j_n \lambda_2^{n-2}  \\
    & \le \sum_{n=2}^{m-1} i_n \lambda_1^{n-2} + \sum_{n=m}^{\f} \lambda_1^{n-2}- \sum_{n=2}^{m-1}i_n \lambda_2^{n-2} \\
    & \le \sum_{n=m}^{\f} \lambda_1^{n-2} = \frac{\lambda_1^{m-2}}{1-\lambda_1},
  \end{align*}
  where the first inequality follows by $i_1\ldots i_{m-1}=j_1\ldots j_{m-1}$, and the second inequality follows by $\la_1<\la_2$.
This, together with $\la_1<\la_2\le \la$,  implies that
 \[ \lambda^{m} \ge \lambda_1^{m-1} \lambda_2 \ge \frac{x(1-\lambda_1 -\lambda_2)}{1-\lambda_2}(\lambda_2 - \lambda_1) \ge x(1-2\lambda)(\lambda_2 - \lambda_1). \]
Therefore,
  \begin{align*}
    |\pi_\lambda(\Psi_x(\lambda_1))- \pi_\lambda(\Psi_x(\lambda_2))| & =  (1-\lambda) \sum_{n=1}^{\f} i_n \lambda^{n-1} - (1-\lambda) \sum_{n=1}^{\f} j_n \lambda^{n-1} \\
    & \ge (1-\lambda)\left( \lambda^{m-1}-\sum_{n=m+1}^{\f} \lambda^{n-1} \right)\\
     &= (1-2\lambda) \lambda^{m-1}  \ge C |\lambda_2 - \lambda_1|,
  \end{align*}
  where $C=x(1-2\lambda)^2/\lambda>0$ (since $\la<1/2$). This proves (\ref{eq:lipschitz-1}), and then completes the proof.
\end{proof}

\begin{lemma}\label{lem:local-dimension-lower-bound}
Let $x\in (0,1/2)$.  If $\lambda \in \Lambda(x) \sm \set{1/2}$ such that $\Psi_x(\lambda)$ does not end with $0^\f$, then for any $\delta >0$,
   \begin{equation*}\label{eq:lower-bound}
    \dim_H (\Lambda(x)\cap [\lambda,\lambda +\delta]) \ge \dim_H K_\lambda.
    \end{equation*}
\end{lemma}
\begin{proof}
Let $\la\in\La(x)\setminus\set{1/2}$ such that $(c_n) = \Psi_x(\lambda)$ does not end with $0^\f$. Take $\delta >0$. We will construct a sequence of subsets in $\La(x)\cap[\la, \la+\de]$ whose Hausdorff dimension can be arbitrarily close to $\dim_H K_\la$.

 Since $\la<1/2$, by Lemma \ref{lem:property-Psi-x} we have  $(c_n)\succ\Psi_x(1/2)=(x_n)$. Then there exists $n_0\ge 2$ such that $c_1\ldots c_{n_0-1}=x_1\ldots x_{n_0-1}$ and  $c_{n_0}>x_{n_0}$. Since $(c_n)$ does not end with $0^\f$, we can find an increasing sequence $\{n_k\}\subset\N$ such that $n_0<n_1<n_2<\cdots$, and $c_{n_k}=1$ for all $k\ge 1$.
  Now for $k \ge 1$, we define
  \begin{equation}\label{eq:lower-bound-(-1)}
  \Omega_{\lambda,k} := \set{ c_1 \ldots c_{n_k-1} 0 i_1i_2\ldots: ~ i_{n+1}\ldots i_{n+k}\ne 0^k~\forall n\ge 0}.
  \end{equation}
  Note by Lemma \ref{lem:property-Psi-x} that $\Psi_x(\Lambda(x) \cap [\lambda, 1/2]) = \set{ (i_n) : (x_n) \preceq (i_n)  \preceq (c_n)}.$
  Then by using $c_{n_0}>x_{n_0}$ and $c_{n_k}=1$ it follows that
 \begin{equation}\label{eq:lower-bound-0}
 \Omega_{\lambda,k} \sse \Psi_x(\Lambda(x) \cap [\lambda, 1/2])\quad\textrm{for all }k\ge 1.
 \end{equation}
Since ${\delta >0}$, by (\ref{eq:lower-bound-(-1)}), (\ref{eq:lower-bound-0}) and  Lemma \ref{lem:property-Psi-x} there exists $N\in\N$ such that
   \begin{equation*}
    \Lambda(x) \cap [\lambda, \lambda+\delta]\supset\Psi_x^{-1}(\Omega_{\lambda,k})\quad \forall k\ge N.
    \end{equation*}
    So, to finish the proof it suffices to prove  that
    \begin{equation}\label{eq:lower-bound-1}
      \lim_{k\to\f}\dim_H\Psi_x^{-1}(\Om_{\la,k})\ge \dim_H K_\la.
    \end{equation}

Take $k\ge N$, and consider the map  $\pi_\lambda \circ \Psi_x: \Psi_x^{-1}(\Omega_{\lambda,k}) \to \pi_\lambda(\Omega_{\lambda,k})$.
  Let $\lambda_1,\lambda_2 \in \Psi_x^{-1}(\Omega_{\lambda,k})$ with $\lambda_1 < \lambda_2$, and write $\Psi_x(\lambda_1) = (i_n),  \Psi_x(\lambda_2) = (j_n).$ Then $(i_n), (j_n)\in\Om_{\la, k}$.
  Since $\la_1<\la_2$, by Lemma \ref{lem:property-Psi-x} we have  $(i_n) \succ (j_n)$. So there exists  $m>n_k$  such that $i_1\ldots i_{m-1}=j_1\ldots j_{m-1}$ and $i_m>j_m$.
  Note   that $i_mi_{m+1}\ldots$ does not contain $k$ consecutive zeros. Then
  \begin{equation}\label{eq:lower-bound-2}
  x=(1-\lambda_1) \sum_{n=1}^{\f} i_n \lambda_1^{n-1}>(1-\lambda_2) \sum_{n=1}^{\f} i_n \lambda_1^{n-1}> (1-\lambda_2)\left( \sum_{n=1}^{m} i_n \lambda_1^{n-1}  + \lambda_1^{m+k-1} \right).
  \end{equation}
On the other hand,
  \begin{equation}\label{eq:lower-bound-3}
  x= (1-\lambda_2) \sum_{n=1}^{\f} j_n \lambda_2^{n-1}\le (1-\lambda_2) \sum_{n=1}^{m} i_n \lambda_2^{n-1}.
  \end{equation}
 Note that $\la_1,\la_2 \in \La(x)\cap[\la, 1/2]$. Then by (\ref{eq:lower-bound-2}) and  (\ref{eq:lower-bound-3})   it follows  that
   \begin{align*}
     \lambda^{m+k-1} \le \lambda_1^{m+k-1} & < \sum_{n=1}^{m} i_n (\lambda_2^{n-1} -\lambda_1^{n-1}) \\
     & < \sum_{n=1}^{\f} (\lambda_2^{n-1} -\lambda_1^{n-1}) = \frac{1}{1-\lambda_2}-\frac{1}{1-\lambda_1}  = \frac{\lambda_2 -\lambda_1}{(1-\lambda_1)(1-\lambda_2)}.
   \end{align*}
This implies that
  \begin{equation}\label{eq:kong-11}
  \begin{split}
    |\pi_\lambda(\Psi_x(\lambda_1))- \pi_\lambda(\Psi_x(\lambda_2))| & =  (1-\lambda) \sum_{n=1}^{\f} i_n \lambda^{n-1} - (1-\lambda) \sum_{n=1}^{\f} j_n \lambda^{n-1} \\
    & \le (1-\lambda)\sum_{n=m}^{\f} \lambda^{n-1} \\
    &= \lambda^{m-1}   < \frac{\lambda_2 -\lambda_1}{\lambda^k (1-\lambda_1)(1-\lambda_2)} {\le} \frac{4}{\lambda^k}(\lambda_2 -\lambda_1),
  \end{split}
  \end{equation}
  where the last inequality follows by $\la_1,\la_2 \le 1/2$.

  So, by (\ref{eq:lower-bound-(-1)}) and (\ref{eq:kong-11}) it follows that
  \begin{align*}
   \dim_H \Psi_x^{-1}(\Omega_{\lambda,k}) &\ge \dim_H \pi_\lambda(\Omega_{\lambda,k}) \\
   & = \dim_H \pi_\lambda\left(\set{(i_n): i_{n+1}\ldots i_{n+k}\ne 0^k~\forall n\ge 0} \right) \\
    & \ge \dim_H  \pi_\lambda \left( \set{(i_n)\in\set{0,1}^\N: i_{n}=1\textrm{ for all }n\equiv 0(\textrm{mod}\; k)} \right) \\
    & = -\frac{(k-1)\log 2}{k \log \lambda}~\; \to \frac{\log 2}{-\log \la}=\dim_H K_\la
  \end{align*}
  as $k\to\f$. This proves (\ref{eq:lower-bound-1}), and then completes the proof.
\end{proof}

\begin{proof}[Proof of Proposition \ref{th:local-dimension}]
Take $\la\in\La(x)$. Note that $\La(x)\subset[x, 1/2]$ and $x, 1/2\in\La(x)$. We will prove (\ref{eq:local-dimension}) in the following two cases.

Case I. $\lambda \in\La(x)\cap[x, 1/2)$. Then by Lemma \ref{lem:local-dim-upper-bound} it follows that for any $\de\in(0, 1/2-\la)$,
\[
\dim_H (\Lambda(x) \cap (\lambda -\delta,\lambda +\delta))\leq \dim_H (\Lambda(x) \cap [x,\la+\delta])\leq \dim _HK_{\lambda +\delta }=\frac{\log 2}{-\log(\la+\de)}.
\]
This implies that
\begin{equation}
  \label{eq:local-dimension-1}
   \lim_{\de\to 0^+}\dim_H (\Lambda(x) \cap (\lambda -\delta,\lambda +\delta))\le \frac{\log 2}{-\log \la}=\dim_H K_\la.
\end{equation}

On the other hand, take $\de >0$. Note by Theorem \ref{main:topology} that $\La(x)$ is a Cantor set, and $\la\in\La(x)$. Then we can find a sequence $ \{\la_k\}$ in $\La(x)\cap(\la-\de, \la+\de)$ such that each $\Psi_x(\la_k)$ does not end with $0^\f$, and $\la_k\to \la$ as $k\to\f$. Therefore, by Lemma \ref{lem:local-dimension-lower-bound} it follows that
\begin{align*}
  \dim_H (\Lambda(x) \cap (\lambda -\delta,\lambda +\delta))&\ge \dim_H(\La(x)\cap[\la_k, \la+\de])\\
  &\ge \dim_H K_{\la_k}=\frac{\log 2}{-\log \la_k}~\; \to\frac{\log 2}{-\log\la}=\dim_H K_\la
\end{align*}
as $k\to\f$.
This, together with (\ref{eq:local-dimension-1}), proves (\ref{eq:local-dimension}).

Case II. $\lambda =1/2$. The proof is similar to that for the second part of Case I.  Let $\de>0$. Since $\La(x)$ is a Cantor set and $\max\La(x)=1/2$, there exists a sequence ${\{\la_k\}}$ in $\La(x)\cap (1/2-\de, 1/2)$ such that each $\Psi_x(\la_k)$ does not end with $0^\f$, and $\la_k\nearrow 1/2$ as $k\to\f$. Then by Lemma \ref{lem:local-dimension-lower-bound} it follows that
\begin{align*}
  \dim_H (\Lambda(x) \cap (1/2 -\delta,1/2 +\delta))&\ge \dim_H(\La(x)\cap[\la_k, \la_{k+1}])\\
  &\ge \dim_H K_{\la_k}=\frac{\log 2}{-\log\la_k}~\; \to 1=\dim_H K_{1/2},
\end{align*}
proving (\ref{eq:local-dimension}).
\end{proof}

 As a direct consequence of Proposition \ref{th:local-dimension} we have the following   result of $\La(x)$.
\begin{corollary}\label{cor:dim-local}
  Let $x \in (0,1/2)$. Then for any open interval $I\subset\R$ with $\Lambda(x) \cap I \ne \emptyset$  we have \[ \dim_H (\Lambda(x) \cap I) = \sup_{\lambda \in \Lambda(x) \cap I} \dim_H K_\lambda. \]
\end{corollary}
\begin{proof}[Proof of Theorem \ref{main:measure-dim}]
By Corollary \ref{cor:dim-local} it follows that
\[
\dim_H \La(x)=\dim_H(\La(x)\cap (x, 1/2))=\sup_{ \la \in \La(x)\cap (x, 1/2)}\dim_H K_\la=1.
\]
Furthermore, for any   $n\in\N$ the Hausdorff dimension of $\La_n(x):=\La(x)\cap[x,1/2-1/n]$ is strictly smaller than one, and thus each $\La_n(x)$ has zero Lebesgue measure. Since $\La(x)\setminus\set{1/2}=\bigcup_{n=1}^\f\La_n(x)$, the set  $\La(x)$ also has zero Lebesgue measure. This together with Proposition \ref{th:local-dimension} completes the proof.
\end{proof}

\section{Hausdorff dimension of the intersection $\bigcap_{i=1}^p\La(y_i)$}\label{sec:intersection}
Given finitely many numbers $y_1, y_2, \ldots, y_p\in(0,1/2)$, we will show in this section that the intersection $\bigcap_{i=1}^p \La(y_i)$ has full Hausdorff dimension (see Theorem \ref{main:intersection}). Note by Theorem \ref{main:topology} that each set $\La(y_i)$ is a Cantor set. We will construct in each $\La(y_i)$ a sequence of Cantor subsets $C_\ell(y_i), \ell\ge 1$ such that each $C_\ell(y_i)$ has the same maximum point $1/2$, and the thickness of $C_\ell(y_i)$ tends to infinity as $\ell\to\f$. Then by using a result from Hunt, Kan and Yorke \cite{Hunt-Kan-Yorke-1992} (see Lemma \ref{theorem-HKY} below) we conclude that the intersection $\bigcap_{i=1}^p \La(y_i)$ contains a sequence of Cantor subsets whose thickness tends to infinity. This, together with Lemma \ref{lemma-thickness-dimension} (see below), implies that $\bigcap_{i=1}^p \La(y_i)$ has full Hausdorff dimension.

\subsection{Thickness {of} a Cantor set in $\R$}
First we recall  the \emph{thickness} of a Cantor set in $\R$ {from} Newhouse \cite{Newhouse-1970} (see \cite{Astels_2000} for some recent progress).
Let $E$ be a Cantor set in $\R$ with its convex hull $conv(E)=E_0$. Then the complement $E_0 \sm E = \bigcup_{n=1}^\f V_n $ is the union of countably many disjoint open intervals.
The sequence $\mscr{V}= (V_1, V_2,\ldots)$ is called a \emph{defining sequence} for $E$.
If moreover $|V_1| \ge |V_2| \ge |V_3| \ge \cdots$, where $|V|$ denotes the diameter of a set $V\subset\R$, then we call $\mscr{V}$ an \emph{ordered defining sequence} for $E$.
Let $E_n := E_0 \sm \bigcup_{k=1}^n V_k$. Then $E_n$ is the union of finitely many disjoint closed intervals. So,
for any $n\ge 1$, the open interval $V_n$ is contained in some connected component of $E_{n-1}$, say $E_{n-1}^*$.
Then the set $E_{n-1}^* \sm V_n$ is the union of two closed intervals $L_{\mscr{V}}(V_n)$ and $R_{\mscr{V}}(V_n)$, where we always assume that $L_{\mscr{V}}(V_n)$ lies to the left of $R_{\mscr{V}}(V_n)$. We emphasize that both intervals $L_{\V}(V_n)$ and $R_{\V}(V_n)$   have positive length, since otherwise $E$ will contain isolated points which is impossible.
Then the {thickness of $E$ with respect to the defining sequence $\mscr{V}$} is  defined by
\begin{equation}\label{eq:def-thickness-def}
 \tau_\mscr{V}(E):= \inf_{n\ge 1} \min\set{ \frac{|L_{\mscr{V}}(V_n)|}{|V_n|}, \frac{|R_{\mscr{V}}(V_n)|}{|V_n|}},
\end{equation}
and the \emph{thickness} of $E$ is defined by
\begin{equation}\label{eq:def-thickness}
 \tau(E) := \sup_{\mscr{V}}\tau_{\mscr{V}}(E),
\end{equation} where the supremum is taken over all defining sequences  $\mscr{V}$ for $E$.
It was shown in \cite{Williams-1991} that $\tau(E)=\tau_\mscr{V}(E)$ for every ordered defining sequence $\mscr{V}$ for $E$.

The following lower bound for the Hausdorff dimension of a Cantor set in $\R$ in terms of thickness was proven by Newhouse \cite{Newhouse-1979} (see also \cite[P.~77]{Palis_Takens_1993}).
\begin{lemma} \cite[P.~107]{Newhouse-1979} \label{lemma-thickness-dimension}
  If $E$ is a Cantor set in $\R$, then  \[ \dim_H E \ge \frac{\log 2}{\log \big(2+\frac{1}{\tau(E)} \big) }. \]
\end{lemma}

Two Cantor sets in $\R$ are called \emph{interleaved} if neither set lies in {the closure of} a gap of the other.
The following result for the intersection of two interleaved Cantor sets was shown by Hunt, Kan and Yorke \cite{Hunt-Kan-Yorke-1992}.
\begin{lemma}\cite[Theorem 1]{Hunt-Kan-Yorke-1992}\label{theorem-HKY}
  There exists a function $\varphi: (1+\sqrt{2} ,\f) \to (0,\f)$ such that for all interleaved Cantor sets $E$ and $F$ {in $\R$} with $\tau(E) , \tau(F) \ge t > 1+ \sqrt{2}$, there exists a Cantor subset $K \sse E \cap F$ with $\tau(K)\ge \varphi(t)$.
\end{lemma}
\begin{remark}\label{rem:intersection}
\begin{enumerate}[{\rm(i)}]
  \item In \cite[P.~882]{Hunt-Kan-Yorke-1992} the authors pointed out that when $t$ is sufficiently large,  $\varphi(t)$  is of order $\sqrt{t}$. So, Lemma \ref{theorem-HKY} implies that if the thicknesses  of two interleaved Cantor sets $E$ and $F$ {in $\R$} are sufficiently large, then the thickness of the resulting Cantor set $K\subset E\cap F$ is also very large.

  \item It is clear that if two Cantor sets $E$ and $F$ {in $\R$} have the same  maximum point $\xi$, then they are interleaved.
  Furthermore, if the maximum point $\xi$ is also an accumulation point of $E\cap F$, then from the proof of   \cite[Theorem 1]{Hunt-Kan-Yorke-1992} (see also \cite[P.~887]{Hunt-Kan-Yorke-1992})  it follows that   the resulting  Cantor set $K\subset E\cap F$ in Lemma \ref{theorem-HKY} can be required to have the same maximum point $\xi$.
\end{enumerate}
\end{remark}

We first construct a sequence of disjoint Cantor subsets of $\Lambda(x)$ whose thickness tends to infinity. Then by the following lemma, we can construct a sequence of Cantor subsets of $\Lambda(x)$ with the same maximum point, whose thickness also tends to infinity.
These Cantor sets will be used to show that the intersection has full Hausdorff dimension.

\begin{lemma}\label{lemma:thickness-infinite}
Let  $\{F_k\}_{k=1}^\f$ be a sequence of Cantor sets with $\alpha_k = \min F_k$ and $\beta_k = \max F_k$. Suppose that
\begin{enumerate}[{\rm(i)}]
\item $\alpha_1 < \beta_1 < \alpha_2 < \beta_2 < \ldots < \alpha_k < \beta_k < \ldots$, and $\beta: = \lim_{k \to \f} \beta_k$;

\item $\lim_{k \to \f} \tau(F_k) =+\f$;

\item $$\lim_{k \to \f} \frac{\beta_k - \alpha_k}{\alpha_{k+1} - \beta_k}=+\f,\quad\lim_{k \to \f} \frac{\beta - \alpha_{k+1}}{\alpha_{k+1} - \beta_k} = +\f.$$
\end{enumerate}
Then  $$\lim_{\ell \to \f} \tau \left( \bigcup_{k=\ell}^\f F_k \cup \{\beta\} \right) = +\f.$$
\end{lemma}
\begin{proof}
For $k \ge 1$, let $\mathscr{V}_k = \{V_{k,j}\}_{j=1}^\f$ be an ordered defining sequence of the Cantor set $F_k$, i.e.,
\begin{equation}\label{eq:defining-sequence-Ek}
  [\alpha_k, \beta_k] \setminus F_k = \bigcup_{j=1}^\f V_{k,j}
\end{equation}
with $|V_{k,1}| \ge |V_{k,2}| \ge |V_{k,3}| \ge \ldots$.
For $\ell \ge 1$, write $$ C_\ell := \bigcup_{k=\ell}^\f F_k \cup \{\beta\}.$$
Note by (i) that each $C_\ell$ is a Cantor set, and the convex hull of $C_\ell$ is $[\alpha_\ell, \beta].$
Moreover, we have $$[\alpha_\ell, \beta] \setminus C_\ell = \bigcup_{k=\ell}^\f (\beta_k, \alpha_{k+1})~\cup~\bigcup_{k=\ell}^\f \bigcup_{j=1}^\f V_{k,j}.$$

To estimate the thickness of $C_\ell$, we enumerate the open intervals $ (\beta_k, \alpha_{k+1})$ and $V_{j,k}$ with $ k \ge \ell, j\ge 1$ in the following way:
\begin{equation}\label{eq:defining-sequence-C-ell}
\begin{gathered}
  \xymatrix@=0.8em{
    (\beta_\ell, \alpha_{\ell+1})\ar[d] & (\beta_{\ell+1}, \alpha_{\ell+2}) \ar"3,1" &  (\beta_{\ell+2}, \alpha_{\ell+3}) \ar"4,1" & (\beta_{\ell+3}, \alpha_{\ell+4}) & \cdots \\
    V_{\ell, 1} \ar[ru] & V_{\ell, 2} \ar[ru] & V_{\ell, 3} \ar[ru] & V_{\ell, 4} & \cdots \\
    V_{\ell+1, 1} \ar[ru] & V_{\ell+1, 2} \ar[ru] & V_{\ell+1, 3} \ar[ru] & V_{\ell+1, 4}& \cdots \\
    V_{\ell+2, 1}\ar[ru] & V_{\ell+2, 2} \ar[ru] & V_{\ell+2, 3} \ar[ru] & V_{\ell+2, 4} & \cdots \\
    \cdots & \cdots & \cdots & \cdots & \cdots
  }
\end{gathered}
\end{equation}
This means that we first remove from $[\al_\ell, 1/2]$ the open interval $(\beta_\ell, \al_{\ell+1})$, and next remove $V_{\ell, 1}$, and then $(\beta_{\ell+1}, \al_{\ell+2}), V_{\ell+1,1}, V_{\ell, 2}, (\beta_{\ell+2}, \al_{\ell+3})$, and so on.
Thus, (\ref{eq:defining-sequence-C-ell}) gives a defining sequence $\mscr{W}_\ell=\set{(\beta_k, \al_{k+1}), V_{k,j}: k\ge \ell, j\ge 1}$ for $C_\ell$ (see Figure 1).

\begin{figure}[h!]
\begin{center}
\begin{tikzpicture}[
    scale=11,
    axis/.style={very thick, ->},
    important line/.style={thick},
    dashed line/.style={dashed, thin},
    pile/.style={thick, ->, >=stealth', shorten <=2pt, shorten
    >=2pt},
    every node/.style={color=black}
    ]

    \draw[important line] (0, 0)node[above,scale=1pt]{$\al_\ell$}--(0.4, 0)node[above,scale=1pt]{$\beta_\ell$};
     \draw[important line] (0, -0.1) --(0.18, -0.1);  \node[above,scale=0.8pt] at(0.23, -0.1){$V_{\ell,1}$}; \draw[important line] (0.28, -0.1) --(0.4, -0.1);
     \draw[important line] (0, -0.2) --(0.08, -0.2);  \node[above,scale=0.7pt] at(0.1, -0.2){$V_{\ell,3}$}; \draw[important line] (0.12, -0.2) --(0.18, -0.2);
     \draw[important line] (0.28, -0.2) --(0.33, -0.2);  \node[above,scale=0.7pt] at(0.345, -0.2){$V_{\ell,2}$}; \draw[important line] (0.36, -0.2) --(0.4, -0.2);

    \node[above,scale=1pt]at(0.2,-0.35){$F_\ell$};

    \draw[important line] (0.55, 0)node[above,scale=1pt]{$\al_{\ell+1}$}--(0.8, 0)node[above,scale=1pt]{$\beta_{\ell+1}$};
    \draw[important line] (0.55, -0.1) --(0.66, -0.1);  \node[above,scale=0.8pt] at(0.68, -0.1){$V_{\ell+1,1}$}; \draw[important line] (0.7, -0.1) --(0.8, -0.1);
     \draw[important line] (0.55, -0.2) --(0.595, -0.2);  \node[above,scale=0.7pt] at(0.61, -0.2){$V_{\ell+1,3}$}; \draw[important line] (0.625, -0.2) --(0.66, -0.2);
     \draw[important line] (0.7, -0.2) --(0.745, -0.2);  \node[above,scale=0.7pt] at(0.755, -0.2){$V_{\ell+1,2}$}; \draw[important line] (0.765, -0.2) --(0.8, -0.2);

      \node[above,scale=1pt]at(0.675,-0.35){$F_{\ell+1}$};

    \draw[important line] (0.9, 0)node[above,scale=1pt]{$\al_{\ell+2}$}--(1.1, 0)node[above,scale=1pt]{$\beta_{\ell+2}$};
    \draw[important line] (0.9, -0.1) --(0.99, -0.1);  \node[above,scale=0.8pt] at(1.01, -0.1){$V_{\ell+2,1}$}; \draw[important line] (1.03, -0.1) --(1.1, -0.1);
     \draw[important line] (0.9, -0.2) --(0.94, -0.2);  \node[above,scale=0.7pt] at(0.95, -0.2){$V_{\ell+2,3}$}; \draw[important line] (0.96, -0.2) --(0.99, -0.2);
     \draw[important line] (1.03, -0.2) --(1.06, -0.2);  \node[above,scale=0.7pt] at(1.07, -0.2){$V_{\ell+2,2}$}; \draw[important line] (1.08, -0.2) --(1.1, -0.2);
 \node[above,scale=1pt]at(1,-0.35){$F_{\ell+2}$};

  \draw[important line, densely dotted] (1.15,0)--(1.3,0)node[above,scale=1pt]{$\beta$};

\end{tikzpicture}
\end{center}
\caption{A defining sequence $\W_\ell=\set{(\beta_{k}, \al_{k+1}), V_{k,j}: k\ge \ell, j\ge 1}$ for the Cantor set $C_\ell=\bigcup_{k=\ell}^\f F_k \cup \{\beta\}$, and for each $k\ge \ell$ a defining sequence $\V_k=\set{V_{k,j}}_{j=1}^\f$ for the Cantor set $F_k$; see (\ref{eq:defining-sequence-Ek}) and (\ref{eq:defining-sequence-C-ell}) for more explanation.}
\label{fig:1}
\end{figure}
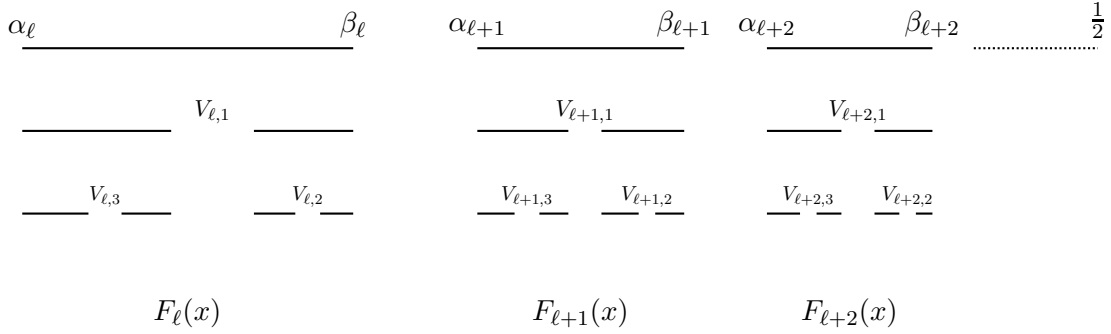

For the defining sequence $\W_\ell$ for $C_\ell$, by (\ref{eq:defining-sequence-C-ell}) we have that
\[ L_{\W_\ell}\big( (\beta_k,\alpha_{k+1}) \big)=[\alpha_k,\beta_k], \quad R_{\W_\ell}\big( (\beta_k,\alpha_{k+1}) \big)= [\alpha_{k+1},\beta]\quad \text{for any } k \ge \ell;\]
and
\[ L_{\W_\ell}(V_{k, j})=L_{\V_k}(V_{k, j}),\quad R_{\W_\ell}(V_{k,j})=R_{\V_k}(V_{k,j}) \quad\text{for any } k\ge \ell, j \ge 1. \]
Note that $\tau_{\V_k}(F_k) = \tau (F_k)$. Then by (\ref{eq:def-thickness-def}) it follows that
\begin{align*}
  \tau_{\W_\ell}(C_\ell) &=\inf_{k\ge \ell, j\ge 1}\min\set{ \frac{\beta_k - \alpha_k}{\alpha_{k+1} - \beta_k}, \frac{\beta - \alpha_{k+1}}{\alpha_{k+1} - \beta_k}, \frac{|L_{\V_k}(V_{k,j})|}{|V_{k,j}|}, \frac{|R_{\V_k}(V_{k,j})|}{|V_{k,j}|}}\\
   &= \inf_{k\ge \ell}\min \set{ \tau(F_k), \frac{\beta_k - \alpha_k}{\alpha_{k+1} - \beta_k}, \frac{\beta - \alpha_{k+1}}{\alpha_{k+1} - \beta_k}}.
\end{align*}
Note that $\tau (C_\ell) \ge \tau_{\W_\ell}(C_\ell)$. We conclude by (ii) and (iii) that $\lim_{\ell \to \f} \tau(C_\ell) = +\f$.
\end{proof}

\subsection{Construction of Cantor subsets of $\La(x)$}
Let $x\in(0,1/2)$.
We will construct a sequence of Cantor subsets $\{F_k(x)\}_{k=1}^\f$ of $\La(x)$ satisfying the assumptions (i)--(iii) in Lemma \ref{lemma:thickness-infinite}.
Note by Lemma \ref{lem:property-Psi-x} that \[\Psi_x(\La(x))=\Om(x)=\set{(i_n): ({x_n})\lle (i_n)\lle 01^\f},\] where $(x_n)=\Psi_x(1/2)$.
Moreover, $\Psi_x$ is a decreasing homeomorphism from $\La(x)$ to $\Om(x)$. Based on this, our strategy to construct these Cantor subsets $\set{F_k(x)}_{k=1}^\f$ in $\La(x)$ is to   construct a sequence of Cantor subsets $\set{\Om_k(x)}_{k=1}^\f$ in the symbolic space $\Om(x)$.

Note by the definition of $(x_n)=\Psi_x(1/2)$ that $x_1=0$ and $(x_n)$ does not end with $1^\f$.
Denote by $\{n_k\}$ the set of all indices $n>1$ such that $x_n=0$.
Then $x_{n_k}=0$ for any $k\ge 1$, and $x_n=1$ for any $n_k<n<n_{k+1}$.
For $k \ge 1$, we define
\begin{equation}\label{eq:cantor-subset-1}
  \Omega_{k}(x) := \set{(i_n): x_1 \ldots x_{n_k -1} 1 0^\f \preceq (i_n) \preceq x_1 \ldots x_{n_k -1} 1^\f }.
\end{equation}
Then ${\Omega_{k}(x)} \sse \Omega(x)$, which implies that
\begin{equation}\label{eq:cantor-subset-2}
  F_k(x) := \Psi_x^{-1}(\Omega_{k}(x))\subset\La(x).
\end{equation}
Note by (\ref{eq:cantor-subset-1}) that for each $k\ge 1$ the set $(\Om_k(x), \rho)$ is a topological Cantor set, where $\rho$ is the metric defined in (\ref{eq:metric-rho}). By Lemma \ref{lem:property-Psi-x} it follows that each $F_k(x)$ is a Cantor subset of $\La(x)$.

Write $\alpha_k = \min F_k(x)$ and $\beta_k = \max F_k(x)$. Clearly, $\al_k$ and $\beta_k$ depend on $x$. For simplicity we will suppress this dependence in our notation if no confusion arises.
Since $\Psi_x$ is decreasing by Lemma \ref{lem:property-Psi-x}, by (\ref{eq:cantor-subset-1}) and (\ref{eq:cantor-subset-2}) we have
\[ \alpha_k = \Psi_x^{-1}(x_1 x_2 \ldots x_{n_k -1} 1^\f),\quad \beta_k = \Psi_x^{-1}(x_1 x_2 \ldots x_{n_k -1} 1 0^\f). \]
Note that $x_1 x_2 \ldots x_{n_{k+1}-1}1^\f=x_1 x_2\ldots x_{n_k-1}01^\f \prec x_1 x_2 \ldots x_{n_k -1} 1 0^\f$ for all $k\ge 1$. Thus, $$\alpha_1 < \beta_1 < \alpha_2 < \beta_2 < \ldots < \alpha_k < \beta_k < \ldots.$$
Furthermore, the sequence $x_1 x_2 \ldots x_{n_k -1} 1 0^\f$ decreases to $(x_n)=\Psi_x(1/2)$ as $k\to\f$, again by Lemma \ref{lem:property-Psi-x} we obtain that $\beta_k \nearrow \frac{1}{2}$ as $k\to\f$.
Therefore, the sequence $\set{F_k(x)}_{k=1}^\f$ of disjoint Cantor sets satisfies the assumption (i) of Lemma  \ref{lemma:thickness-infinite}.

In the following  we show that the sequence $\set{F_k(x)}_{k=1}^\f$ also satisfies the assumption (ii) of Lemma \ref{lemma:thickness-infinite}.
\begin{proposition}\label{proposition:thickness-F-k}
 For any $x\in(0,1/2)$ we have
 $\lim_{k\to\f}\tau(F_k(x))=+\f.$
\end{proposition}
In view of the definition of thickness,  we first describe a defining sequence for the Cantor set $F_k(x)$. Clearly, $[\alpha_k,\beta_k]$ is the convex hull of $F_k(x)$.
By Lemma \ref{lem:property-Psi-x} it follows that
\[ [\alpha_k ,\beta_k ] \sm F_k(x) = \bigcup_{\omega \in \set{0,1}^*} V_{k, \om}, \]
where
\[ V_{k, \om}  :=\left( \Psi_x^{-1}(x_1 x_2 \ldots x_{n_k -1} 1\ \omega\ 10^\f) , \Psi_x^{-1}(x_1 x_2 \ldots x_{n_k -1} 1\ \omega\ 0 1^\f) \right). \]
We enumerate these open intervals $V_{k, \om}$, $\omega \in \set{0,1}^*$ according first to the length of $\om$ and then to the lexicographical order of $\om$:
\begin{equation}\label{eq-5-4}
  \begin{split}
    &V_{k,\epsilon};\quad V_{k, 0}, V_{k, 1};\quad  V_{k,00}, V_{k,01}, V_{k,10}, V_{k,11}; \\
    &V_{k,000}, V_{k, 001}, V_{k,010}, V_{k, 011}, V_{k,100}, V_{k, 101}, V_{k, 110}, V_{k, 111};~~ \ldots,
  \end{split}
\end{equation}
where $\epsilon$ is the empty word.
Thus, (\ref{eq-5-4}) gives a defining sequence $\V_{k} = \{ V_{k, \om}: \omega \in \set{0,1}^* \}$ for $F_k(x)$, and moreover, we have
\begin{equation}\label{eq:F-k-L-R}
\begin{split}
  L_{\V_k}(V_{k,\om}) & = \left[ \Psi_x^{-1}(x_1 x_2 \ldots x_{n_k -1} 1\ \omega\ 1^\f) , \Psi_x^{-1}(x_1 x_2 \ldots x_{n_k -1} 1\ \omega\ 1 0^\f) \right],  \\
  R_{\V_k}(V_{k,\om}) & =\left[ \Psi_x^{-1}(x_1 x_2 \ldots x_{n_k -1} 1\ \omega\ 0 1^\f) , \Psi_x^{-1}(x_1 x_2 \ldots x_{n_k -1} 1\ \omega\ 0^\f) \right].
\end{split}
\end{equation}

Note by (\ref{eq:F-k-L-R}) that the endpoints of $L_{\V_k}(V_{k,\om})$ and $ R_{\V_k}(V_{k,\om})$ have codings ending with $\omega1^\f, \omega0^\f, \omega10^\f, \omega01^\f$, respectively. In view of (\ref{eq:def-thickness-def}), to prove $\lim_{k\to\f}\tau(F_k(x))$ we need the following inequalities.

\begin{lemma}\label{lemma-case-1}
Let $x\in(0,1/2)\setminus\set{1/4}$ with $(x_n)=\Psi_x(1/2)$, and let $m\ge 3$ such that $x_m=1$.
\begin{enumerate}[{\rm(i)}]
\item  If $\la_1, \la_2\in\La(x)$ satisfy
$\Psi_x(\la_1)= j_1 j_2 \ldots j_q 1^\f$ and $\Psi_x(\lambda_2) =   j_1 j_2 \ldots j_q 0^\f$,
  then\[ \lambda_2 - \lambda_1 \ge \frac{1}{4} \lambda_2^q. \]

\item  If $\la_3, \la_4\in\La(x)$ satisfy $\Psi_x(\la_3)=x_1 \ldots x_m j_1 \ldots j_q 1 0^\f$ and $\Psi_x(\la_4)=x_1 \ldots x_m j_1  \ldots j_q 0 1^\f $,
  then
   \[ \lambda_4 - \lambda_3 \le \min\set{2(1-2\lambda_3) \lambda_3^{q+2},  2(1-2\lambda_4) \frac{\lambda_4^{m+q}}{\lambda_3^{m-2}}}. \]
  \end{enumerate}
\end{lemma}
Before giving the proof we emphasize that for $x=1/4$ we have $(x_n)=\Psi_x(1/2)=010^\f$. So we can not find $x_m=1$ for $m\ge 3$, which plays an essential role in the proof of Lemma \ref{lemma-case-1} (ii).
\begin{proof}
Since $x\in(0,1/2)\setminus\set{1/4}$, the expansion $(x_n)=\Psi_x(1/2)$ satisfies $x_m=1$ for some $m\ge 3$.
 For (i), let $\la_1, \la_2\in\La(x)$ with $\Psi_x(\la_1)= j_1 j_2 \ldots j_q 1^\f$ and $\Psi_x(\lambda_2) = j_1 j_2 \ldots j_q 0^\f$. Then by Lemma \ref{lem:property-Psi-x} we have $\la_1<\la_2$.
 Note that  $x = \pi_{\lambda_1}( j_1 j_2 \ldots j_q 1^\f) = \pi_{\lambda_2}( j_1 j_2 \ldots j_q 0^\f)$. Then
  \begin{align*}
    \lambda_2^{q}  = \pi_{\lambda_2}(0^{q}1^\f)
     & = \pi_{\lambda_2}( j_1 j_2 \ldots j_q 1^\f) - \pi_{\lambda_2}(j_1 j_2 \ldots j_q  0^\f) \\
     & =  \pi_{\lambda_2}( j_1 j_2 \ldots j_q 1^\f) - \pi_{\lambda_1}( j_1 j_2 \ldots j_q 1^\f) \\
     & \le 4 (\lambda_2 - \lambda_1),
  \end{align*}
  where the inequality follows by Lemma \ref{lem:property-Pi} (ii) since for any $(i_n)\in\Om(x)$ we have
  \[ \frac{\mrm{d} \pi_\lambda((i_n))}{\mrm{d} \lambda}= \sum_{n=2}^{\f} n \left( \frac{n-1}{n} -\lambda \right)  i_n\lambda^{n-2}\le \sum_{n=2}^{\f} \frac{n-1}{2^{n-2}}=4. \]
  This proves (i).

 For (ii) let $\la_3, \la_4\in\La(x)$ such that
 \[\Psi_x(\la_3)=x_1 \ldots x_m j_1 \ldots j_q 1 0^\f\quad\textrm{and}\quad \Psi_x(\la_4)=x_1 \ldots x_m j_1  \ldots j_q 0 1^\f,\] where $x_m=1$  with $m\ge 3$.
 Then by Lemma \ref{lem:property-Psi-x} we have $\lambda_3 < \lambda_4$. Note that $x = \pi_{\lambda_3}(x_1 \ldots x_{m} j_1 j_2 \ldots j_q 1 0^\f) = \pi_{\lambda_4}(x_1 \ldots x_{m} j_1 j_2 \ldots j_q 0 1^\f)$.
  Then
\begin{equation}\label{eq:estimate-0}
  \begin{split}
    &~\quad \pi_{\lambda_3}(0^{m+q} 1 0^\f) - \pi_{\lambda_4}(0^{m+q+1} 1^\f) \\
    &=  \pi_{\lambda_4}(x_1 \ldots x_{m} j_1 j_2 \ldots j_q 0^\f) - \pi_{\lambda_3}(x_1 \ldots x_{m} j_1 j_2 \ldots j_q 0^\f) \\
    &\ge  \frac{1}{2} \lambda_3^{m-2} (\lambda_4-\lambda_3),
  \end{split}
\end{equation}
where the last inequality follows by Lemma \ref{lem:property-Pi} (ii) since by using $x_m=1$ with $m\ge 3$ and (\ref{deriv}) we have
\[ \frac{\mrm{d} \pi_\lambda(x_1 \ldots x_{m} j_1 j_2 \ldots j_q 0^\f)}{\mrm{d} \lambda}\ge m \left( \frac{m-1}{m} -\lambda \right) x_m \lambda^{m-2} \ge \frac{1}{2} \lambda^{m-2}. \]
Therefore, by (\ref{eq:estimate-0}) and using $\la_3<\la_4$ it follows that
  \begin{align*}
    \lambda_4 - \lambda_3 & \le \frac{2}{\lambda_3^{m-2}} \left( (1-\lambda_3)\lambda_3^{m+q} - \lambda_4^{m +q +1} \right) \\
    & \le \frac{2}{\lambda_3^{m-2}} \left( (1-\lambda_3)\lambda_3^{m+q} - \lambda_3^{m +q +1} \right)
     = 2(1-2\lambda_3) \lambda_3^{q+2},
  \end{align*}
and
  \begin{align*}
    \lambda_4- \lambda_3 & \le \frac{2}{\lambda_3^{m-2}} \left( (1-\lambda_3)\lambda_3^{m+q} - \lambda_4^{m +q +1} \right) \\
    & \le \frac{2}{\lambda_3^{m-2}} \left( (1-\lambda_4)\lambda_4^{m+q} - \lambda_4^{m +q +1} \right)
     = 2(1-2\lambda_4) \frac{\lambda_4^{m+q}}{\lambda_3^{m-2}}.
  \end{align*}
This completes the proof.
\end{proof}

When $x=1/4$ we prove similar inequalities by using  different estimation.
\begin{lemma}\label{lemma-case-2}
Let $x=1/4$. Then $(x_n)=\Psi_x(1/2)=010^\f$.
\begin{enumerate}[{\rm(i)}]
  \item  If $\la_1, \la_2\in\La(x)$ satisfy $\Psi_x(\la_1)=010^m j_1 j_2 \ldots j_q 1^\f$ and $\Psi_x(\la_2)=010^m j_1 j_2 \ldots j_q 0^\f$,
  then  \[ \lambda_2 - \lambda_1 \ge \frac{ \lambda_2^{m+2+q} }{1-2\lambda_1 + (m+3) 2^{-m}}. \]
  \item  If $\la_3, \la_4\in\La(x)$ satisfy $\Psi_x(\la_3)= 01 j_1 j_2 \ldots j_q 1 0^\f$ and $\Psi_x(\la_4)=01 j_1 j_2 \ldots j_q 0 1^\f$,
  then \[ \lambda_4 - \lambda_3 \le \lambda_3^{2+q}. \]
  \end{enumerate}
\end{lemma}

\begin{proof}
For (i) we note by Lemma \ref{lem:property-Psi-x} that $\lambda_1 < \lambda_2$. Since  $x = \pi_{\lambda_1}( 010^m j_1 j_2 \ldots j_q 1^\f) = \pi_{\lambda_2}( 010^m j_1 j_2 \ldots j_q 0^\f),$
  we have
  \begin{align*}
    \lambda_2^{m+2+q}  = \pi_{\lambda_2}(0^{m+2+q}1^\f)
     &= \pi_{\lambda_2}(010^m j_1 j_2 \ldots j_q 1^\f) - \pi_{\lambda_2}(010^m j_1 j_2 \ldots j_q  0^\f) \\
     & =  \pi_{\lambda_2}( 010^m j_1 j_2 \ldots j_q 1^\f) - \pi_{\lambda_1}( 010^m j_1 j_2 \ldots j_q 1^\f) \\
     & \le \big( 1-2\lambda_1 + (m+3) 2^{-m} \big) (\lambda_2 - \lambda_1),
  \end{align*}
 where the inequality follows by Lemma \ref{lem:property-Pi} (ii) since by (\ref{deriv}) we have
 \begin{align*}
    \frac{\mrm{d} \pi_\lambda(010^m j_1 j_2 \ldots j_q 1^\f)}{\mrm{d} \lambda}
    & \le (1-2\lambda)+\sum_{n=m+3}^{\f} n \left( \frac{n-1}{n} -\lambda \right) \lambda^{n-2} \\
    & \le 1-2\lambda+\sum_{n=m+3}^{\f} \frac{n-1}{2^{n-2}}  = 1- 2\lambda + (m+3) 2^{-m}.
 \end{align*}
This proves (i).

For (ii), note by Lemma \ref{lem:property-Psi-x} that $\lambda_3 < \lambda_4$. Then by using $ x = \pi_{\lambda_3}(01 j_1 j_2 \ldots j_q 1 0^\f) = \pi_{\lambda_4}(01 j_1 j_2 \ldots j_q 0 1^\f)$ it follows that
   we have
  \begin{align*}
     \pi_{\lambda_3}(0^{2+q} 1 0^\f) - \pi_{\lambda_4}(0^{3+q} 1^\f)
    = & \pi_{\lambda_4}( 01 j_1 j_2 \ldots j_q 0^\f) - \pi_{\lambda_3}(01 j_1 j_2 \ldots j_q 0^\f) \\
    \ge & \pi_{\lambda_4}( 01  0^\f) - \pi_{\lambda_3}(01  0^\f) \\
    =&   (1-\lambda_3 -\lambda_4)(\lambda_4- \lambda_3).
  \end{align*}
This implies that
  \begin{align*}
    \lambda_4 - \lambda_3 & \le \frac{1}{1-\lambda_3 -\lambda_4} \left( (1-\lambda_3)\lambda_3^{2+q} - \lambda_4^{3 +q} \right) \\
    & \le \frac{1}{1-\lambda_3 -\lambda_4} \left( (1-\lambda_3)\lambda_3^{2+q} - \lambda_4\lambda_3^{2 +q} \right)  = \lambda_3^{2+q},
  \end{align*}
  as desired.
\end{proof}

\begin{proof}[Proof of Proposition \ref{proposition:thickness-F-k}]
Write $\Psi_x(1/2)=(x_n)$. Let $\{n_k\}$ be the enumeration of all indices $n>1$ such that $x_n=0$.
Recall from (\ref{eq-5-4}) the defining sequence $\V_k = \{V_{k,\om}: \om \in \{0,1\}^*\}$ for $F_k(x)$.
It follows from (\ref{eq:def-thickness-def}) that $$\tau_{\V_k}(F_k(x)) = \inf\left\{ \frac{|L_{\V_k}(V_{k,\om})|}{|V_{k,\om}|}, \frac{|R_{\V_k}(V_{k,\om})|}{|V_{k,\om}|}: \om \in \{0,1\}^* \right\}.$$
Note by (\ref{eq:F-k-L-R}) that
\begin{align*}
  L_{\V_k}(V_{k,\om}) & = \left[ \Psi_x^{-1}(x_1 x_2 \ldots x_{n_k -1} 1\ \omega\ 1^\f) , \Psi_x^{-1}(x_1 x_2 \ldots x_{n_k -1} 1\ \omega\ 1 0^\f) \right] =: \big[ \ga_{\om,1}, \ga_{\om,2} \big], \\
  R_{\V_k}(V_{k,\om}) & =\left[ \Psi_x^{-1}(x_1 x_2 \ldots x_{n_k -1} 1\ \omega\ 0 1^\f) , \Psi_x^{-1}(x_1 x_2 \ldots x_{n_k -1} 1\ \omega\ 0^\f) \right] =: \big[ \ga_{\om,3},\ga_{\om,4} \big].
\end{align*}
Then \[ V_{k,\om} =\left( \Psi_x^{-1}(x_1 x_2 \ldots x_{n_k -1} 1\ \omega\ 10^\f) , \Psi_x^{-1}(x_1 x_2 \ldots x_{n_k -1} 1\ \omega\ 0 1^\f) \right)=\big( \ga_{\om,2},\ga_{\om,3} \big).\]
Observe by Lemma \ref{lem:property-Psi-x} that $\al_k \le \ga_{\om,1} < \ga_{\om,2} < \ga_{\om,3} < \ga_{\om,4}$, where $\al_k=\min F_k(x)$.

\medskip
\noindent
{\bf Case (A). $x\in(0,1/2)\setminus\{1/4\}$.}
Then there exists an integer $m \ge 3$ such that $x_m=1$.
Take a sufficiently large $k$ so that $n_k > m$.
By Lemma \ref{lemma-case-1} it follows that
\[ \ga_{\om,2} - \ga_{\om,1} \ge \frac{1}{4} \ga_{\om,2}^{n_k + q +1}, \quad \ga_{\om,4} - \ga_{\om,3} \ge \frac{1}{4} \ga_{\om,4}^{n_k +q +1} ,\]
and
\[ \ga_{\om,3} - \ga_{\om,2} \le 2(1-2\ga_{\om,2}) \ga_{\om,2}^{n_k +q -m +2}, \]
where $q$ is the length of the word $\om$.
So, for any $\om \in \{0,1\}^*$ we obtain that
\begin{align*}
\frac{|L_{\V_k}(V_{k,\om})|}{|V_{k,\om}|} &= \frac{\ga_{\om,2} - \ga_{\om,1}}{\ga_{\om,3} - \ga_{\om,2}} \ge \frac{\ga_{\om,2}^{m-1}}{8(1-2\ga_{\om,2})} \ge \frac{\alpha_k^{m-1}}{8(1-2\alpha_k)}, \\
 \frac{|R_{\V_k}(V_{k,\om})|}{|V_{k,\om}|} &= \frac{\ga_{\om,4} - \ga_{\om,3}}{\ga_{\om,3} - \ga_{\om,2}} \ge \frac{\ga_{\om,4}^{m-1}}{8(1-2\ga_{\om,2})} \ge \frac{\alpha_k^{m-1}}{8(1-2\alpha_k)}.
\end{align*}
Thus, for sufficiently large $k$ we have
\[ \tau(F_k(x)) \ge \tau_{\mscr{V}_k}(F_k(x)) \ge \frac{\alpha_k^{m-1}}{8(1-2\alpha_k)}. \]
Note that $\alpha_k\nearrow 1/2$ as $k \to \f$ because $\beta_k \nearrow 1/2$.
Therefore,  $\tau(F_k(x)) \to +\f$ as $k \to \f$.

\medskip
{\bf \noindent Case (B). $x= 1/4$.} Then $(x_n)=\Phi_x(1/2) =  010^\f$, which gives
  $x_1\ldots x_{n_k-1}=010^{n_k-3}$ for all $k \ge 1$.
By Lemma \ref{lemma-case-2} it follows that
\[ \ga_{\om,2} - \ga_{\om,1} \ge \frac{ \ga_{\om,2}^{n_k+q+1} }{1-2\ga_{\om,1} + n_k 2^{3-n_k}}, \quad \ga_{\om,4} - \ga_{\om,3} \ge \frac{ \ga_{\om,4}^{n_k+q+1} }{1-2\ga_{\om,3} + n_k 2^{3-n_k}} ,\]
and \[ \ga_{\om,3} - \ga_{\om,2} \le \ga_{\om,2}^{n_k +q}, \]
where $q$ is the length of the word $\om$.
This implies that for all $\om \in \{0,1\}^*$,
\begin{align*}
 \frac{|L_{\V_k}(V_{k,\om})|}{|V_{k,\om}|} &= \frac{\ga_{\om,2} - \ga_{\om,1}}{\ga_{\om,3} - \ga_{\om,2}} \ge \frac{\ga_{\om,2}}{1-2\ga_{\om,1} + n_k 2^{3-n_k}} \ge \frac{\alpha_k}{1-2\alpha_k + n_k 2^{3-n_k}}, \\
 \frac{|R_{\V_k}(V_{k,\om})|}{|V_{k,\om}|} &= \frac{\ga_{\om,4} - \ga_{\om,3}}{\ga_{\om,3} - \ga_{\om,2}} \ge \frac{\ga_{\om,4}}{1-2\ga_{\om,3} + n_k 2^{3-n_k}} \ge \frac{\alpha_k}{1-2\alpha_k + n_k 2^{3-n_k}}.
\end{align*}
Thus,
\[ \tau(F_k(x)) \ge \tau_{\mscr{V}_k}(F_k(x)) \ge \frac{\alpha_k}{1-2\alpha_k + n_k 2^{3-n_k}}. \]
Note that $\alpha_k \nearrow 1/2$ and $n_k \to +\f$ as $k \to \f$.
Therefore, $\tau(F_k(x)) \to +\f$ as $k \to \f$, completing the proof.
\end{proof}

\subsection{Hausdorff dimension of $\bigcap_{i=1}^p \La(y_i)$}
For any $\ell\ge 1$, we define
\begin{equation}\label{eq:cantor-subset}
 C_\ell(x) := \bigcup_{k=\ell}^\f F_k(x) \cup \{1/2\}.
\end{equation}
Note that $\mathrm{conv}(F_k(x))=[\al_k,\beta_k]$ for all $k\ge 1$, and
\[\al_1<\beta_1<\al_2<\beta_2<\cdots<\al_k<\beta_k<\cdots,\quad\textrm{and}\quad \lim_{k\to\f}\beta_k= 1/2. \]
Since $F_k(x)\subset\La(x)$ for all $k\ge 1$ and $1/2\in\La(x)$, it follows that each $C_\ell(x)$ is a Cantor subset of $\La(x)$.

\begin{proposition}\label{th:thickness-C-ell}
 For any  $x\in(0,1/2)$ we have
 $\lim_{\ell\to\f}\tau(C_\ell(x))=+\f.$
\end{proposition}
\begin{proof}
Take $x\in(0,1/2)$. Write $\Psi_x(1/2)=(x_n)$. Let $\{n_k\}$ be the enumeration of all indices $n>1$ such that $x_n=0$.
By Lemma \ref{lemma:thickness-infinite} and Proposition \ref{proposition:thickness-F-k}, it remains to prove that
\begin{equation}\label{eq:remain-estimation}
  \lim_{k \to \f} \frac{\beta_k - \alpha_k}{\alpha_{k+1} - \beta_k} = + \f \quad\text{and}\quad \lim_{k \to \f} \frac{1/2 - \alpha_{k+1}}{\alpha_{k+1} - \beta_k} = + \f,
\end{equation}
where
\begin{equation}\label{eq:k1}
 \alpha_k = \Psi_x^{-1}(x_1 x_2 \ldots x_{n_k -1} 1^\f),\quad \beta_k = \Psi_x^{-1}(x_1 x_2 \ldots x_{n_k -1} 1 0^\f),
\end{equation}
and
\begin{equation}\label{eq:k2}
  \alpha_{k+1} =\Psi_x^{-1}(x_1 \ldots x_{n_k -1} x_{n_k} \ldots x_{n_{k+1}-1}1^\f )= \Psi_x^{-1}(x_1 x_2 \ldots x_{n_k -1} 0 1^\f).
\end{equation}

\medskip

{\bf\noindent Case (A). $x\in(0,1/2)\setminus\set{1/4}$.}
There exists $m\ge 3$ such that $x_m=1$.
For sufficiently large $k$, we have $n_k > m$.
By (\ref{eq:k1}), (\ref{eq:k2}) and Lemma \ref{lemma-case-1} it follows that
\[ \beta_k - \alpha_k \ge \frac{1}{4} \beta_k^{n_k},\quad  \alpha_{k+1} - \beta_k \le  \min\left\{ 2(1-2\beta_k) \beta_k^{n_k -m +1}, 2(1-2\alpha_{k+1}) \frac{\alpha_{k+1}^{n_k-1 }}{\beta_k^{m-2}} \right\}, \]
which implies
\[ \frac{\beta_k - \alpha_k}{\alpha_{k+1} - \beta_k} \ge \frac{\beta_k^{m-1}}{8(1-2\beta_k)} \quad\text{and}\quad \frac{1/2 - \alpha_{k+1}}{\alpha_{k+1} - \beta_k} \ge \frac{ \beta_k^{m-2} }{4\alpha_{k+1}^{n_k -1} }. \]
Note that $\alpha_k \nearrow 1/2$, $\beta_k \nearrow 1/2$ and $n_k \to \f$ as $k \to \f$.
We obtain (\ref{eq:remain-estimation}) by letting $k \to \f$.

\medskip
{\bf \noindent Case (B). $x= 1/4$.} Then $(x_n)=\Phi_x(1/2) =  010^\f$.
By (\ref{eq:k1}), (\ref{eq:k2}) and Lemma \ref{lemma-case-2} it follows that
\begin{equation}\label{eq:kong}
 \beta_k - \alpha_k \ge \frac{ \beta_k^{n_k} }{1-2\alpha_k + n_k 2^{3-n_k}}, \quad \alpha_{k+1} - \beta_k \le  \beta_k^{n_k  -1}.
\end{equation}
This implies
\begin{equation}\label{eq:kong-2}
 \frac{\beta_k - \alpha_k}{\alpha_{k+1} - \beta_k} \ge \frac{\beta_k}{1-2\alpha_k + n_k 2^{3-n_k}}.
\end{equation}

On the other hand, note by (\ref{eq:k2}) that \[ \frac{1}{4} = x = \pi_{\alpha_{k+1}} (x_1 x_2 \ldots x_{n_k -1} 0 1^\f) =\pi_{\al_{k+1}}(010^{n_k-2}1^\f)= (1-\alpha_{k+1})\alpha_{k+1} + \alpha_{k+1}^{n_k}. \]
This implies that
\begin{equation}\label{eq:kong'}
\left( \frac{1}{2} -\alpha_{k+1} \right)^2 = \alpha_{k+1}^{n_k}, \quad \text{ and thus}\quad \frac{1}{2} -\alpha_{k+1} = \alpha_{k+1}^{n_k/2}.
\end{equation}
Note by (\ref{eq:kong}) that $\alpha_{k+1} - \beta_k \le \beta_k^{n_k -1} { <} \alpha_{k+1}^{n_k-1}$.
This, together with (\ref{eq:kong'}), implies that
\begin{equation} \label{eq:kong-3}
\frac{1/2 - \alpha_{k+1}}{\alpha_{k+1} - \beta_k} \ge \frac{ 1 }{\alpha_{k+1}^{n_k/2 -1} } .
\end{equation}
Note that $\alpha_k\nearrow 1/2$, $\beta_k\nearrow 1/2$ and $n_k \to +\f$ as $k\to\f$.
Letting $k\to\f$ in (\ref{eq:kong-2}) and (\ref{eq:kong-3}) we obtain (\ref{eq:remain-estimation}), completing the proof.
\end{proof}

\begin{proof}[Proof of Theorem \ref{main:intersection}]
 Let $y_1, \ldots, y_p\in(0,1/2)$. Then by (\ref{eq:def-thickness}), (\ref{eq:cantor-subset}) and Proposition \ref{th:thickness-C-ell} it follows that each $\La(y_i)$ contains a sequence of Cantor subsets $C_\ell(y_i), \ell\ge 1$ such that $\max C_\ell(y_i)=1/2$ for all $\ell\ge 1$, and the thickness $\tau(C_\ell(y_i))\to +\f$ as $\ell\to\f$.
 So, by Lemma \ref{theorem-HKY} and Remark \ref{rem:intersection} (i) it follows that for sufficiently large $\ell$ and for any $i, j\in\set{1,2, \ldots, p}$ the intersection $C_\ell(y_i)\cap C_\ell(y_j)$ contains a Cantor subset $C_\ell(y_i, y_j)$ such that
 \begin{equation}\label{eq:kjh-1}
 \tau(C_\ell(y_i, y_j))\to +\f\quad \textrm{as }\ell\to\f.
 \end{equation}
 Note that for any $k\in\set{1,2,\ldots,  p}$ we have $\min C_\ell(y_k)\nearrow 1/2=\max C_\ell(y_k)$  as $\ell\to\f$. {Furthermore, $C_{\ell}(y_k)\supset C_{\ell+1}(y_k)$ for any $\ell\ge 1$.} Then by (\ref{eq:kjh-1}) it follows that the maximum point $1/2$ is an accumulation point of $C_\ell(y_i)\cap C_\ell(y_j)$. So, by Remark \ref{rem:intersection} (ii)  we can require that the resulting Cantor set $C_\ell(y_i, y_j)\subset C_\ell(y_i)\cap C_\ell(y_j)$ has the maximum point $1/2$ for sufficiently large $\ell$ and any $i, j\in\set{1,2,\ldots, p}$.

 Proceeding this argument for all $y_1, y_2,\ldots, y_p$ we obtain that for sufficiently large $\ell$ the intersection $\bigcap_{i=1}^p C_\ell(y_i)$ contains a Cantor subset $C_\ell(y_1,\ldots, y_p)$ such that $\max C_\ell(y_1,\ldots, y_p)=1/2$, and the thickness $\tau(C_\ell(y_1, \ldots, y_p))\to +\f$ as $\ell\to\f$. Therefore, by Lemma \ref{lemma-thickness-dimension} it follows that
 \[
 \dim_H\bigcap_{i=1}^p\La(y_i)\ge\dim_H\bigcap_{i=1}^p C_\ell(y_i)\ge\dim_H C_\ell(y_1,\ldots, y_p)\ge\frac{\log 2}{\log(2+\frac{1}{\tau(C_\ell(y_1, \ldots, y_p))})}\to 1,
 \]
 as $\ell\to\f$. This completes the proof.
\end{proof}

At the end of this section we remark that in the proof of Theorem \ref{main:intersection} we construct a sequence of Cantor subsets $C_\ell(y_1,\ldots, y_p), \ell\ge 1$ in the intersection $\bigcap_{i=1}^p\La(y_i)$ such that the thickness $\tau(C_\ell(y_1,\ldots, y_p))\to +\f$ as $\ell\to\f$. By a recent work of Yavicoli \cite[remark of Theorem 4]{Yavicoli-2020} it follows that the intersection $\bigcap_{i=1}^p\La(y_i)$ contains arbitrarily long arithmetic progression.

\section{Final remarks}\label{sec:final-remarks}
At the end of this paper we point out that our  results Theorem \ref{main:topology}--\ref{main:intersection} can be extended to higher dimensions. To  illustrate this we give two examples.
\begin{example}\label{ex:1}
For $\lambda \in (0, 1/2]$ let $K_\la$ be the self-similar set defined in (\ref{eq:self-similar-set}). Then for $n\in\N$  the product set $\bigotimes_{i=1}^n K_\lambda$ is also a self-similar set in $\R^n$.
For $\mathbf a=(a_1,\ldots, a_n) \in (0,1/2)^n$ let
$$
 \Lambda (\mathbf a):=\{\lambda \in (0, 1/2]: \mathbf a\in \bigotimes_{i=1}^n K_\lambda \}
$$
be the set of parameters $\la\in(0,1/2]$ such that the $n$-dimensional self-similar set $\bigotimes_{i=1}^n K_\la$ contains the given point $\mathbf a$. It is clear that $\La(\mathbf a)=\bigcap_{i=1}^n\Lambda (a_i)$,
where {for $x\in(0,1/2)$ the set} $\La(x)$ is defined as in (\ref{eq:definition-La}). So, by Theorems \ref{main:measure-dim} and \ref{main:intersection} it follows that
$\La(\mathbf a)$ has zero Lebesgue measure and full Hausdorff dimension for any $\mathbf a\in(0,1/2)^n$.
\end{example}

\begin{example}
  \label{ex:2}
  Let $n\in\N$ and let $\la_i\in(0,1/2]$ for all  $1\le i\le n$. Then  $\bigotimes_{i=1}^n K_{\lambda _i} \subset\R^n$ is a  self-affine set generated by the IFS $\set{(\la_1 x_1, \la_2 x_2,\ldots, \la_n x_n)+{\bf i}:\; {\bf i}\in \bigotimes_{i=1}^n\set{0, 1-\la_i} }.$
For any $\mathbf b=(b_1,\ldots, b_n)\in (0,1/2)^n$ let
$$
\La'(\mathbf b):= \{(\lambda _1, \lambda _2,\ldots, \la_n)\in (0, 1/2]^n: \mathbf b\in \bigotimes_{i=1}^n K_{\lambda _i}\}.
$$
Then $(\la_1,\la_2,\ldots, \la_n)\in \La'(\mathbf b)$ if and only if $b_i\in K_{\la_i}$ for all $1\le i\le n$, which is also equivalent to $\la_i\in\La(b_i)$ for all $1\le i\le n$.
So, $\La'(\mathbf b)=\bigotimes_{i=1}^n \Lambda (b_i).$ By Theorem \ref{main:topology} it follows that for any $\mathbf b\in(0,1/2)^n$ the set $\La'(\mathbf b)$ is a Cantor set in $\R^n$, i.e., it is a non-empty compact, totally disconnected and perfect set in $\R^n$.
Furthermore, by \cite[product formula 7.2]{Falconer_1990} and Theorem \ref{main:measure-dim} we obtain that $\La'(\mathbf b)$ has Lebesgue measure zero and $\dim _H\Lambda'(\mathbf b)=n$ for any $\mathbf b\in(0,1/2)^n$.
\end{example}

Note that in Examples \ref{ex:1} and \ref{ex:2} the higher dimensional self-similar sets all have the product form $\bigotimes_{i=1}^n K_{\la_i}$. However, if the higher dimensional self-similar sets   do not have the product form,  our results Theorem  \ref{main:topology}--\ref{main:intersection} can not be applied directly. Some recent progress on the extension of Newhouse thickness theorems to higher dimensions may be useful in this direction (cf.~\cite{Biebler-2020,Falconer-Yavicoli-2022,Feng-Wu-2021,Yavicoli-2023}).

\section*{Acknowledgements}
The authors thank the anonymous referee for many useful comments which greatly improve the paper.
K.~Jiang was supported by NSFC No.~11701302,  Zhejiang Provincial NSF No.~LY20A010009, and the K.C. Wong Magna Fund in Ningbo University.
{D.~Kong was supported by NSFC No.~11971079.}
W.~Li was supported by NSFC No.~12071148 and Science and Technology Commission of Shanghai Municipality (STCSM) No.~22DZ2229014.
Z. Wang was supported by Fundamental Research Funds for the Central Universities No.~YBNLTS2023-016.

\end{document}